\documentclass[12pt]{amsart}
\usepackage{amssymb}
\usepackage{amsmath}
\usepackage{mathrsfs}
\usepackage{versions}
\usepackage{amsfonts}
\usepackage[usenames]{color}
\usepackage{graphicx}
 \newtheorem{theorem}{Theorem}[section]
 
 \newtheorem{lemma}[theorem]{Lemma}

\newtheorem{definition}[theorem]{Definition}
\newtheorem{remark}[theorem]{Remark}

\newtheorem{fact*}{Fact}

\DeclareMathOperator{\diag}{diag}

\DeclareMathOperator\spa{span}

\newcommand\half{{\tfrac 12}}

\newcommand\dd{\mathrm d}

\newcommand\ess{\mathscr{S}}

\newcommand{\T}{\mathbb{T}}

\newcommand{\D}{\mathbb{D}}
\newcommand{\C}{\mathbb{C}}

\newcommand{\norm}[1]{\left\Vert#1\right\Vert}

\newcommand{\ip}[2]{\left\langle #1, #2 \right\rangle}

\newcommand{\inv}{^{-1}}

\newcommand{\ph}{\varphi}

\newcommand\ga{\gamma}

\newcommand\la{\lambda}

\newcommand\si{\sigma}
\newcommand\beq{\begin{equation}}

\newcommand\eeq{\end{equation}}
\newcommand\df{\stackrel{\rm def}{=}}

\newcommand\bbm{\begin{bmatrix}}
\newcommand\ebm{\end{bmatrix}}
\newcommand\bpm{\begin{pmatrix}}
\newcommand\epm{\end{pmatrix}}
\numberwithin{equation}{section}

\newcommand\m{\mathcal{M}}

\newcommand\calk{\mathcal{K}}
\newcommand\call{\mathcal{L}}
\newcommand\calm{\mathcal{M}}
\newcommand\calh{\mathcal{H}}

\newcommand\be{\begin{equation}}
\newcommand\ee{\end{equation}}

\begin{document}
\title{Realization of functions on the symmetrized bidisc}
\author{Jim Agler}
\address{Department of Mathematics, University of California at San Diego, CA \textup{92103}, USA}
\thanks{Partially supported by National Science Foundation grant
DMS 1361720, Engineering and Physical Sciences Research Council grant EP/N03242X/1 and by the London Mathematical Society grant 41527}

\author{N. J. Young}
\address{School of Mathematics and Statistics, Newcastle University, Newcastle upon Tyne NE1 7RU, U.K.
{\em and} School of Mathematics, Leeds University,  Leeds LS2 9JT, U.K.}
\email{Nicholas.Young@ncl.ac.uk}
\date{2nd April, 2017}
\begin{abstract} 
We prove a realization formula and a model formula for analytic functions with modulus bounded by $1$ on the symmetrized bidisc
\[
G\df \{(z+w,zw): |z|<1, \, |w| < 1\}.
\]
As an application we prove a Pick-type theorem giving a criterion for the existence of such a function satisfying a finite set of interpolation conditions.
\end{abstract}
\maketitle

\section{Introduction} \label{intro}

The fascination of the symmetrized bidisc $G$
lies in the fact that much of the classical function theory of the disc $\D$ and bidisc $\D^2$ generalizes in an explicit way
to $G$, but with some surprising twists.  The original motivation for the study of $G$ was its connection
with the spectral Nevanlinna-Pick problem \cite{AY,AY2001}, wherefore the emphasis was on analytic maps from
the unit disc $\D$ into $G$.  However, in studying such maps one is inevitably drawn into studying maps
from $G$ to $\D$; indeed, the duality between these two classes of maps is a central feature of the
theory of hyperbolic complex spaces in the sense of Kobayashi \cite{koba}.  

The idea of a {\em realization formula} for a class of functions has proved potent in both engineering
and operator theory.  Out of hundreds of papers on this topic in the mathematical literature alone, we mention  \cite{NF,Helton1974,Helton,ag1990,BallBolot,BallHuaman,BMV,DM}.  The simplest realization formula provides an elegant connection between
function theory (the Schur class of the disc) and contractive operators on Hilbert space.  It is as follows.

Let $f$ be an analytic function on $\D$ such that $|f(z)|\leq 1$ for all $z\in\D$.  There exists a Hilbert space $\calm$,
a scalar $A\in\C$, vectors $\beta,\ga \in\calm$ and an operator $D$ on $\calm$ such that the operator
\be\label{contr}
\bbm A & 1\otimes\beta\\ \ga\otimes 1& D \ebm \qquad \mbox{ is a contraction on } \C \oplus \calm
\ee
and, for all $z\in\D$,
\be\label{classicform}
f(z)= A  + \ip{z(1-Dz)\inv \ga}{\beta}_\calm.
\ee
Conversely, any function $f$ on $\D$ expressible in the form \eqref{contr}, \eqref{classicform} is an analytic function in
$\D$ satisfying $|f|\leq 1$ on $\D$.

In an earlier paper \cite{AYY} we
gave a realization formula for analytic maps from $\D$ to the closure of  $G$; in this paper we present
the dual notion, a realization formula for analytic maps from $G$ to $\D^-$.

For any open set $U \subset \C^d$ the set of analytic functions on $U$ with values in the closed unit disc $\D^-$
is called the {\em Schur class of $U$} and is denoted by $\ess(U)$.  

We shall use superscripts to denote the components of points in $\C^d$.

For any point $s=(s^1,s^2)\in G$ and any contractive linear operator $T$ on a Hilbert space $\calm$, we define the operator
\beq\label{defsU}
s_T=(2s^2T-s^1)(2-s^1T)\inv \quad \mbox{ on } \calm.
\eeq
Note that $|s^1|<2$ for $s\in G$, and therefore the inverse in equation \eqref{defsU} exists.

We shall derive both `model formulae' and a realization formula for functions in $\ess(G)$.  The latter is the following.

\begin{theorem}\label{main}
Let $\ph\in\ess(G)$.  There exist a Hilbert space $\calm$ and unitary operators
\beq\label{TABCD}
T \mbox{ on }\calm \quad \mbox{ and } \quad \bbm A&B\\C&D \ebm \mbox{ on } \C\oplus \calm
\eeq
such that, for all $s\in G$,
\beq\label{realform1}
\ph(s)=A+Bs_T(1-Ds_T)\inv C.
\eeq

Conversely, any function $\ph$ on $G$ expressible by the formula \eqref{realform1}, where $T, A, B, C,D$ are such that the operators in formula \eqref{TABCD} are unitary, is an analytic function from $G$ to $\D^-$.
\end{theorem}
Both instances of the word `unitary' in the above theorem can validly be replaced by `contractive'.

The classical realization formula \eqref{classicform} is in terms of a single unitary operator (or contraction), whereas our formula for functions in $\ess(G)$ requires the {\em pair} of unitaries (or contractions) \eqref{TABCD}; this is a consequence of the fact that our derivation  invokes two separate lurking isometry arguments. 

The model formula for functions in $\ess(G)$ is derived in Section \ref{model} from the known model formula for $\ess(\D^2)$ by a symmetrization argument.   The realization formula is then deduced from the model formula in Section \ref{theformula}.  A second model formula, involving an integral with respect to a spectral measure, is proved in Section \ref{secondmodel}.  Finally a Pick-type interpolation theorem, giving a solvability criterion for interpolation problems in $\ess(G)$, is demonstrated in Section \ref{pickthm}.  We also give a realization  
formula for bounded analytic {\em operator}-valued functions on $G$.  The proof requires only notational changes from that of Theorem \ref{main}.

This paper is based on a short course of lectures \cite{ag2014} given by the first-named author at the International Centre for the Mathematical Sciences in Edinburgh in 2014.

Two sources for basic facts about the function theory and geometry of $G$ are  \cite[Chapter 7]{jp} and \cite[Appendix A]{aly2016}.

Many authors have generalized the classical realization formula \eqref{contr} to bounded functions on domains other than the disc.  The paper \cite{ag1990} first made it clear that the appropriate class of holomorphic functions for realization theory on certain more general domains $\Omega$ is a subclass of $\ess(\Omega)$, which has become known as the {\em Schur-Agler class} of $\Omega$.  For the disc, the bidisc and the symmetrized bidisc the Schur-Agler class coincides with the Schur class, and so we have no need for its definition in this paper.

The fact that the Schur and Schur-Agler classes of $G$ are equal was proved in \cite{AY} (see also \cite{AY2000}) with the aid of Ando's Theorem on commuting pairs of contractions and a symmetrization argument.  In this paper we show that essentially the same argument, only with a different ending, yields a realization formula for functions in $\ess(G)$.  We believe that the symmetrization argument is a significant item in the toolkit of realization theory.

Model formulae and realization formulae for the Schur-Agler class of $\Omega$, for any domain $\Omega$ having a matrix-polynomial or even a holomorphic operator-valued defining function, are given in \cite{BallBolot,BMV}, together with several applications.  
The question therefore arises as to whether $G$ has a holomorphic operator-valued defining function, and accordingly whether a realization formula for the Schur-Agler class of $G$ can be simply deduced from a known general result.  More specifically, is there a continuous operator-valued function $F$ on the closure $G^-$ of $G$, holomorphic in $G$, which defines $G$ in the following sense?
\beq\label{holdef}
G=\{s\in G^-: \|F(s)\| < 1\}.
\eeq
If so one immediately obtains a realization formula for the general function $\ph$ in the Schur-Agler class of $G$ of the form
\[
\ph(s)=A+ BF(s)(1-DF(s))\inv C
\]
for some contractive (or unitary) operator colligation $ABCD$.  It is therefore significant for this paper that the symmetrized bidisc cannot be defined by a {\em matrix}-valued holomorphic function \cite{KosY}, nor is it known to be defined by an operator-valued holomorphic function.    
We say a little more about this question at the end of the paper.

A generalization of the realization theory of the polydisc to much more general domains, based on {\em test functions}, has been developed by Dritschel, McCullough and others \cite{DM,DMM,BallHuaman}.  We thank a referee for the observation that a realization formula for functions in the Schur-Agler class of $G$ can be derived from the `abstract realization theorem'  \cite[Theorem 2.2]{DM} by the choice of the functions
\[
s\mapsto \frac {2\la s^2-s^1}{2-\la s^1}
\]
(for $|\la|<1$)
as the test functions on $G$.   This procedure is essentially carried out in \cite{BS}, where a 
 realization formula somewhat similar to ours is given \cite[Realization theorem, page 5]{BS}.  However, this approach only yields a realization formula for the Schur-Agler class, not the Schur class, and so to prove Theorem \ref{main} in this way one must invoke \cite{AY}, implicitly utilizing the symmetrization argument we use in this paper.

We are grateful to an anonymous referee for some very helpful remarks which helped us to improve the presentation of this paper.

\section{A model formula for $G$}\label{model}
The notion of a Hilbert space model for a function on the polydisc was introduced in \cite{ag1990}.  A {\em model on} $\D^2$ is a pair $(\calm, u)$ where $\calm=(\calm^1,\calm^2)$ is a pair of Hilbert spaces and $u=(u^1,u^2)$ is a pair of analytic maps from $\D^2$ to $\calm^1, \calm^2$ respectively.  If  $\ph$ is a function on $\D^2$ then $(\calm,u)$ is a {\em model of} $\ph$ if, for all $\la,\mu\in\D^2$,
\beq\label{modelbidisc}
1-\overline{\ph(\mu)}\ph(\la)=(1-\overline{\mu^1}\la^1)\ip{u^1(\la)}{u^1(\mu)}_{\calm^1} + (1-\overline{\mu^2}\la^2)\ip{u^2(\la)}{u^2(\mu)}_{\calm^2}.
\eeq

It is shown in \cite{ag1990} that a function $\ph$ on $\D^2$ belongs to the Schur class $\ess(\D^2)$ if and only if $\ph$ has a model.
In this section we shall adapt the notion of model to $G$ and prove an analogous result by means of a symmetrization argument.

\begin{definition}\label{defGmodel}
A $G$-\emph{model} for a function $\ph$ on $G$ is a triple $(\m,T,u)$ where $\calm$ is a Hilbert space, $T$ is a contraction acting on $\calm$ and $u:G \to \m$ is an analytic function such that, for all $s,t\in G$,
\beq\label{modelform}
 1-\overline{\ph(t)}\ph(s)= \ip{ (1-t_T^* s_T) u(s)}{u(t)}_\calm.
\eeq
\end{definition}
The following is the main result of this section.
\begin{theorem}\label{modelGthm}
Let $\ph$ be a function on $G$.  The following three statements are equivalent.
\begin{enumerate}
\item   $\ph\in\ess(G)$;
\item $\ph$ has a $G$-model;
\item $\ph$ has a $G$-model $(\calm, T, u)$ in which $T$ is a unitary operator on $\calm$.
\end{enumerate}
\end{theorem}
\begin{proof} (2)$\Rightarrow$(1).
  Suppose $\ph$ has a $G$-model $(\calm,T,u)$.  By holding $t$ fixed in equation \eqref{modelform} one can deduce that $\ph$ is analytic on $G$, and on choosing $t=s$ one has
\beq\label{1minus}
1-|\ph(s)|^2=\ip{(1-s_T^*s_T)u(s)}{u(s)}.
\eeq
Now for $s\in G$ we have $|s^1|<2$ and so the function
\[
f_s(\la)=\frac{2\la s^2-s^1}{2-\la s^1}
\]
is analytic for $\la$ in a neighborhood of $\D^-$. Moreover $|f_s|$ is bounded by $1$ on $\D$ \cite[Theorem 2.1, (1)$ \Rightarrow $(4)]{AY04}.  
By von Neumann's inequality $f_s(T)$ is a contraction, that is, $\|s_T\|\leq 1$.  Hence, by equation \eqref{1minus}, $|\ph(s)|\leq 1$.

(3)$\Rightarrow$(2) is trivial.
To prove that (1)$\Rightarrow$(3) we first symmetrize the model \eqref{modelbidisc} for the Schur class of the bidisc.
Denote by superscript $\si$ the transposition of co-ordinates in $\C^2$, so that
\[
(\la^1,\la^2)^\si=(\la^2,\la^1).
\]
Say that a function $h$ on $\D^2\times\D^2$ is {\em doubly symmetric} if it is symmetric with respect to $\si$ in each variable separately, that is, if
\[
h(\la,\mu)=h(\la^\si,\mu)=h(\la,\mu^\si)
\]
for all $\la,\mu\in\D^2$.

A doubly symmetric function $h(\la,\mu)$ on $\D^2\times\D^2$ that is analytic in $\la$ and $\bar\mu$ can be written in terms of the elementary symmetric functions $\la^1+\la^2, \ \la^1\la^2, \ \mu^1+\mu^2$ and $\mu^1\mu^2$.  Specifically,  if $h$ has a Hilbert space model on the bidisc in the sense of the next proposition, then it induces a function on $G$ having a Hilbert space model of the following form.

\begin{lemma} \label{dblysymm}
Let $h$ be a doubly symmetric function on $\D^2 \times \D^2$ such that there exists a model $(\calh,u)$ on the bidisc satisfying, for all
$\la,\mu\in\D^2$,
\beq\label{4.10}
h(\la,\mu)=(1-\overline{\mu^1}\la^1)\ip{u^1(\la)}{u^1(\mu)}_{\calh^1} + (1-\overline{\mu^2}\la^2)\ip{u^2(\la)}{u^2(\mu)}_{\calh^2}.
\eeq
Then there exist a Hilbert space $\calm$, a unitary operator $U$ on $\calm$ and an analytic function $x:G \to \calm$ satisfying
\be\label{hmodel}
h(\lambda,{\mu}) =\Big\langle \left[1 - \bar{t^2}s^2 -\tfrac 12 (\bar{t^1} - s^1\bar{t^2})U -\tfrac12({s^1} - \bar{t^1}{s^2}) U^*\right] x(s),x(t)\Big\rangle_\calm
\ee
for all $\lambda,\mu \in \D^2$, where $s=(\lambda^1+\lambda^2,\lambda^1\lambda^2)$ and $t=(\mu^1+\mu^2,\mu^1\mu^2)$.
\end{lemma}

\begin{proof}
We shall write $u_\la$ in place of $u(\la)$ throughout the proof.

Replace $\lambda$ with $\lambda^\sigma$ and $\mu$ with $\mu^\sigma$ in equation \eqref{4.10} to deduce that
\be\label{4.20}
h(\lambda,{\mu}) =(1-\bar{\mu^2}\lambda^2)\ip{u^1_{\lambda^\sigma}}{u^1_{\mu^\sigma}}
+(1-\bar{\mu^1}\lambda^1)\ip{u^2_{\lambda^\sigma}}{u^2_{\mu^\sigma}}
\ee
for all $\la,\mu\in\D^2$.
On averaging equations \eqref{4.10} and \eqref{4.20} we obtain
\be\label{4.40}
h(\lambda,{\mu}) =\half(1-\bar{\mu^1}\lambda^1)
\ip{\begin{bmatrix}u^1_\lambda\\u^2_{\lambda^\sigma}\end{bmatrix}}
{\begin{bmatrix}u^1_\mu\\u^2_{\mu^\sigma}\end{bmatrix}}
+\half(1-\bar{\mu^2}\lambda^2)
\ip{\begin{bmatrix}u^1_{\lambda^\sigma}\\u^2_\lambda\end{bmatrix}}
{\begin{bmatrix}u^1_{\mu^\sigma}\\u^2_\mu\end{bmatrix}}.
\ee
For every $\la\in\D^2$ define $v_\la\in\calh^1\oplus\calh^2$ by
\[
v_\lambda = \begin{bmatrix}u^1_\lambda\\u^2_{\lambda^\sigma}\end{bmatrix}
\]
for $\la\in\D^2$.   Equation \eqref{4.40} becomes
\be\label{4.60}
\ h(\lambda,{\mu}) =\half(1-\bar{\mu^1}\lambda^1)
\ip{v_\lambda}{v_\mu}
+\half(1-\bar{\mu^2}\lambda^2)
\ip{v_{\lambda^\sigma}}{v_{\mu^\sigma}}.
\ee
 So far we have only used the `weak symmetry' $h(\la^\si,\mu^\si)=h(\la,\mu)$.  Now use the hypothesis $h(\lambda^\sigma,{\mu}) =h(\lambda,\mu)$. On substituting into equation {\eqref{4.60}} we deduce that
\[
(1-\bar{\mu^1}\lambda^1)
\ip{v_\lambda}{v_\mu}
+(1-\bar{\mu^2}\lambda^2)
\ip{v_{\lambda^\sigma}}{v_{\mu^\sigma}}=(1-\bar{\mu^1}\lambda^2)
\ip{v_{\lambda^\sigma}}{v_\mu}
+(1-\bar{\mu^2}\lambda^1)
\ip{v_{\lambda}}{v_{\mu^\sigma}}.
\]
Rearrange the terms in this formula to obtain
\begin{align*}
\ip{v_\lambda}{v_\mu}& + \ip{v_{\lambda^\sigma}}{v_{\mu^\sigma}}
-\ip{v_{\lambda^\sigma}}{v_\mu}-\ip{v_{\lambda}}{v_{\mu^\sigma}} =\\
	&\bar{\mu^1}\lambda^1
\ip{v_\lambda}{v_\mu}
+\bar{\mu^2}\lambda^2
\ip{v_{\lambda^\sigma}}{v_{\mu^\sigma}}-\bar{\mu^1}\lambda^2
\ip{v_{\lambda^\sigma}}{v_\mu}
-\bar{\mu^2}\lambda^1
\ip{v_{\lambda}}{v_{\mu^\sigma}}.
\end{align*}
Both sides of this equation factor, to yield
\be\label{4.80}
\ip{v_\lambda-v_{\lambda^\sigma}}{v_\mu-v_{\mu^\sigma}} = \ip{\lambda^1v_\lambda-\lambda^2v_{\lambda^\sigma}}{\mu^1v_\mu-\mu^2v_{\mu^\sigma}}.
\ee
In other words, the Gramian in $\calh^1\oplus\calh^2$ of the family of vectors $\{v_\la-v_{\la^\si}:\la\in\D^2\}$ is equal to the Gramian of the family $\{\lambda^1v_\lambda-\lambda^2v_{\lambda^\sigma}:\la\in\D^2\}$.   Hence there exists a linear isometry 
\[
L: \overline{\spa}\{v_\la-v_{\la^\si}:\la\in\D^2\} \to \overline{\spa}\{\la^1v_\la-\la^2v_{\la^\si}:\la\in\D^2\}
\]
 such that
\be\label{4.90}
L(v_\lambda-v_{\lambda^\sigma}) = \lambda^1v_\lambda-\lambda^2v_{\lambda^\sigma}
\ee
for all $\lambda \in \D^2$.  Extend $L$ to a unitary operator $U$ on a Hilbert space $\calm \supseteq \calh^1\oplus\calh^2$.

 Rearrange equation \eqref{4.90} (with $L$ replaced by $U$) to obtain
\[
(U-\lambda^1) v_\lambda = (U-\lambda^2)v_{\lambda^\sigma}
\]
or equivalently,
\beq\label{symv}
(U-\lambda^2)^{-1} v_\lambda\ = (U-\lambda^1)^{-1}v_{\lambda^\sigma}.
\eeq
Therefore,  if we define $w_\lambda$ by the formula
\beq\label{defw}
w_\lambda = (U-\lambda^2)^{-1} v_\lambda
\eeq
then
\be\label{4.91}
v_\lambda = (U-\lambda^2) w_\lambda\ \ \text{ and }\ \ v_{\lambda^\sigma} = (U-\lambda^1)w_{\lambda}.
\ee
If we substitute these formulae  into equation \eqref{4.60} we obtain 
\begin{align}
 h(\lambda,{\mu}) &=\half(1-\bar{\mu^1}\lambda^1)
\ip{(U-\lambda^2) w_\lambda}{(U-\mu^2) w_\mu} \notag \\
	& \hspace*{1in} +\half(1-\bar{\mu^2}\lambda^2) \ip{(U-\lambda^1)w_{\lambda}}{(U-\mu^1)w_\mu}\notag\\ 
	&=\half(1-\bar{\mu^1}\lambda^1)
\ip{(U-\mu^2)^*(U-\lambda^2) w_\lambda}{ w_\mu}\notag \\
	&\hspace*{1in} +\half(1-\bar{\mu^2}\lambda^2)
\ip{(U-\mu^1)^*(U-\lambda^1)w_{\lambda}}{w_\mu}\notag\\
	&= \ip{Z w_\la}{w_\mu} \label{rformh}
\end{align}
where
\beq\label{4.95}
{Z} = \half \Big[(1-\bar{\mu^1}\lambda^1)(U-\mu^2)^*(U-\lambda^2) +(1-\bar{\mu^2}\lambda^2)(U-\mu^1)^*(U-\lambda^1)\Big].
\eeq

Gathering terms in equation \eqref{4.95} we find that
\begin{align*}
{Z}=(1 -& \bar{\mu^1}\bar{\mu^2}\lambda^1\lambda^2) -\half\Big(\overline{\mu^1 +\mu^2} - (\lambda^1+\lambda^2)\overline{\mu^1\mu^2}\Big)\ U\\
	&- \half\Big(\lambda^1 +\lambda^2 - \overline{ \mu^1+\mu^2}{\lambda^1\lambda^2}\Big)\ U^*
\end{align*}
which, in the symmetric variables
\beq\label{svar}
s^1 = \lambda^1+\lambda^2, \qquad  s^2 = \lambda^1 \lambda^2
\eeq
and
\beq\label{tvar}
t^1=\mu^1+\mu^2 , \qquad  t^2 = \mu^1 \mu^2
\eeq
becomes
\be\label{4.96}
{Z}=1 - \bar{t^2}s^2-\half  \big(\bar{t^1} - s^1\bar{t^2}\big)\ U - \half \big({s^1} - \bar{t^1}{s^2}\big)\ U^*.
\ee
Hence
\beq\label{4.97}
h(\la,\mu)= \ip{(1 - \bar{t^2}s^2-\half  \big(\bar{t^1} - s^1\bar{t^2}\big)\ U - \half \big({s^1} - \bar{t^1}{s^2}\big)\ U^*)w_\la}{w_\mu}
\eeq
for all $\la,\mu\in\D^2$.

From the definition \eqref{defw} of $w_\la$ it is clear that $w:\D^2\to \calm$ is analytic, and from equation \eqref{symv} we have
\[
w_{\la^\si}=(U-\la^1)\inv v_{\la^\si}= (U-\la^2)\inv v_\la = w_\la.
\]
Thus $w$, being symmetric, factors through $G$: there exists an analytic function $x:G\to\calm$ such that, for all $\la\in\D^2$,
\[
w_\la=x(\la^1+\la^2,\la^1\la^2)=x(s^1,s^2).
\]

On combining this equation with equation \eqref{4.97} we obtain the desired model formula \eqref{hmodel} for $h$.
\end{proof}

We resume the proof of (1)$\Rightarrow$(3) in Theorem \ref{modelGthm}.
Let $\ph\in\ess(G)$.  The function
\[
\tilde\ph(\la)=\ph(\la^1+\la^2,\la^1\la^2)
\]
belongs to $\ess(\D^2)$, and therefore, by \cite[Theorem 1.12]{ag1990}, has a model $(\calh,v)$ on $\D^2$, which is to say that
\[
1-\overline{\tilde\ph(\mu)}\tilde\ph(\la)=(1-\overline{\mu^1}\la^1)\ip{v^1(\la)}{v^1(\mu)}_{\calh^1} + (1-\overline{\mu^2}\la^2)\ip{v^2(\la)}{v^2(\mu)}_{\calh^2}
\]
for all $\la,\mu \in \D^2$.  The left hand side of this equation is clearly a doubly symmetric function of $(\la,\mu)$, and so, by Lemma \ref{dblysymm}, there exist a Hilbert space $\calm$, a unitary operator $U$ on $\calm$ and an analytic function $x:G \to \calm$ satisfying (in terms of the variables $s,t$ defined in equations \eqref{svar} and \eqref{tvar})
\be\label{phmodel}
1-\overline{\ph(t)}\ph(s) =\Big\langle \left[1 - \bar{t^2}s^2 -\tfrac 12 (\bar{t^1} - s^1\bar{t^2})U -\tfrac12({s^1} - \bar{t^1}{s^2}) U^*\right] x(s),x(t)\Big\rangle_\calm.
\ee
By inspection,
\[
1 - \bar{t^2}s^2 -\tfrac 12 (\bar{t^1} - s^1\bar{t^2})U -\tfrac12(s^1 - \bar{t^1}{s^2}) U^*=
(1-\tfrac12 t^1U)^*(1-\tfrac12 s^1U) - (t^2U-\tfrac12 t^1)^*(s^2U-\tfrac12 s^1).
\]
In the notation $s_U$ introduced in Definition \ref{defsU},
\[
s_U=(s^2U-\tfrac12 s^1)(1-\tfrac12 s^1U)\inv
\]
and we have
\[
1 - \bar{t^2}s^2 -\tfrac 12 (\bar{t^1} - s^1\bar{t^2})U -\tfrac12(s^1 - \bar{t^1}{s^2}) U^*=
(1-\tfrac12 t^1U)^*(1-t_U^*s_U)(1-\tfrac12 s^1U) .
\]

For $s\in G$ let 
\[
u(s)=(1-\half s^1U)x(s).
\]
Then $u:G\to\calm$ is analytic, and  equation \eqref{phmodel} can be written
\begin{align*}
1-\overline{\ph(t)}\ph(s) &=\ip{  (1-\tfrac12 t^1U)^*(1-t_U^*s_U)(1-\tfrac12 s^1U)  x(s)}{x(t)}\\
	&=\ip{(1-t_U^*s_U)u(s)}{u(t)}.
\end{align*}
Thus $(\calm, U, u)$ is a $G$-model for $\ph$.  Therefore (1)$\Rightarrow$(3).
\end{proof}

There is an analogue of Theorem \ref{modelGthm} for operator-valued functions. It is proved by making only notational
changes in the above proof.  If $\calh,\calk$ are Hilbert spaces, $\call(\calh,\calk)$ is the Banach space of bounded linear operators from $\calh$ to $\calk$ in the operator norm, then we define the corresponding Schur class $\ess(G;\calh,\calk)$ to be the set of analytic maps $\ph: G\to \call(\calh,\calk)$ such that $\ph(\la)$ is a contraction for all $\la\in G$.  The notion of $G$-model is extended as follows.
\begin{definition}\label{defGmodelop}
A $G$-\emph{model} for an operator-valued function $\ph:G \to  \call(\calh,\calk)$ is a triple $(\m,T,u)$ where $\calm$ is a Hilbert space, $T$ is a contraction acting on $\calm$ and $u:G \to \call(\calh,\calm)$ is an analytic function such that, for all $s,t\in G$,
\beq\label{modelformop}
 1-\ph(t)^*\ph(s)= u(t)^*(1-t_T^* s_T) u(s).
\eeq
\end{definition}
The generalization of Theorem \ref{modelGthm} is then:
\begin{theorem}\label{modelGthmop}
Let $\ph$ be a function from $G$ to $\calh,\calk)$.  The following three statements are equivalent.
\begin{enumerate}
\item   $\ph\in\ess(G;\calh,\calk)$;
\item $\ph$ has a $G$-model;
\item $\ph$ has a $G$-model $(\calm, T, u)$ in which $T$ is a unitary operator on $\calm$.
\end{enumerate}
\end{theorem}
Details of the proof of this theorem can be found in \cite[Lemmas 3.3 and 3.4]{AY}, where the result was used to derive a certain integral representation formula \cite[Theorem 3.5]{AY} and thereafter to show that if $G^-$ is a spectral set for a commuting pair of operators then $G^-$ is a {\em complete} spectral set.  The fact that the Schur and Schur-Agler classes of $G$ coincide then follows by standard manoeuvres based on the Arveson Extension and Stinespring Representation Theorems.
In this paper we use  Theorem \ref{modelGthm} and its analogue for operator-valued functions to take a more direct route to realization formulae for $\ess(G)$ and $\ess(G;\calh,\calk)$ (Theorems \ref{therealznthm} and \ref{therealznthmop} in the next section).

\section{The realization formula}\label{theformula}
There is a standard way to deduce a realization formula from a model formula with the aid of a `lurking isometry' argument.
We shall apply such an argument to derive the following slight strengthening of Theorem \ref{main} in the introduction.
\begin{theorem}\label{therealznthm}
Let $\ph\in\ess(G)$.  There exist a scalar $A$, a Hilbert space $\calm$, vectors $\beta,\gamma \in\calm$ and operators $D, U$ on $\calm$ such that $U$ is unitary, the operator
\beq\label{collig}
\bbm A & 1\otimes\beta \\ \ga\otimes 1& D \ebm \mbox{ is unitary on } \C\oplus\calm
\eeq
and, for all $s\in G$,
\beq\label{thisisit}
\ph(s)=A+\ip{s_U(1-Ds_U)\inv \ga}{\beta}_\calm.
\eeq

Conversely, if a scalar $A$, a Hilbert space $\calm$, vectors $\beta,\ga \in\calm$ and operators $T, D$ on $\calm$ are given such that $T$ is a contraction and 
\beq\label{colligbis}
\bbm A & 1\otimes\beta \\ \ga\otimes 1& D \ebm \mbox{ is a contraction on } \C\oplus\calm
\eeq
 then the function $\ph$ on $G$ defined by
\beq\label{lfsT}
\ph(s)=A+\ip{s_T(1-Ds_T)\inv\ga}{\beta}
\eeq
belongs to $\ess(G)$.
\end{theorem}
\begin{proof}
Let $\ph\in\ess(G)$.  By Theorem \ref{modelGthm}, $\ph$ has a $G$-model $(\calm, U,u)$ where $U$ is a unitary operator on $\calm$.  By the definition of a $G$-model we have
\[
1-\overline{\ph(t)}\ph(s) =\ip{(1-t_U^*s_U)u(s)}{u(t)}
\]
for all $s,t\in G$.  Rearrange to obtain
\[
1+\ip{s_Uu(s)}{t_Uu(t)}=\overline{\ph(t)}\ph(s)+\ip{u(s)}{u(t)},
\]
which is to say that the two families of vectors
\[
\bpm 1 \\ s_Uu(s)\epm_{s\in G} \quad \mbox{ and }\quad \bpm \ph(s) \\ u(s) \epm_{s\in G}
\]
in $\C\oplus\calm$ have the same Gramians.  Hence there exists an isometry
\beq\label{propL}
L:\overline{\spa} \left\{ \bpm 1 \\ s_Uu(s)\epm_{s\in G}\right\}  \to \overline{\spa}\left\{ \bpm \ph(s) \\ u(s) \epm_{s\in G}\right\}
\eeq
 such that
\[
L\bpm 1 \\ s_Uu(s)\epm = \bpm \ph(s) \\ u(s) \epm
\]
for every $s\in G$.  If necessary enlarge the Hilbert space $\calm$ (and simultaneously the unitary operator $U$ on $\calm$) so that the isometry $L$ extends to a unitary operator
\[
 L^\sharp \sim \bbm A & 1\otimes\beta \\ \ga\otimes 1& D \ebm \mbox{ on } \C\oplus\calm
\]
for some vectors $\beta,\ga \in\calm$.  By equation \eqref{propL}, for any $s\in G$,
\begin{align}\label{2eq}
A+\ip{s_U u(s)}{\beta}&=\ph(s), \notag\\
\ga+Ds_U u(s)&=u(s).
\end{align}
Now $s_U=f_s(U)$ where 
\[
f_s(\la)=\frac {2\la s^2-s^1}{2-\la s^1}
\]
 for $\la $ in a neighborhood of $\D^-$.   The linear fractional map $f_s$ maps $\D$ onto the open disc with centre and radius
\[
 2\frac{\overline{s^1}s^2-s^1}{4-|s^1|^2} \quad \mbox{ and } \quad \frac{|(s^1)^2-4s^2|}{4-|s^1|^2}.
\]
Therefore, by von Neumann's inequality,
\[
\|s_U\| \leq \sup_\D |f_s| = \frac{2|s^1-\bar{s^1}s^2|+|(s^1)^2-4s^2|}{4-|s^1|^2}.
\]  
But, by \cite[Theorem 2.1]{AY04}, the right hand side of this equation is less than one for $s\in G$.
Hence $1-Ds_U$ is invertible for any $s\in G$, and we may eliminate $u(s)$ from equations \eqref{2eq} to obtain the realization formula \eqref{thisisit} for $\ph(s)$.

The converse statement is easy since equation \eqref{lfsT} expresses $\ph(s)$ as a linear fractional transform of the contraction $s_T$ with a contractive coefficient matrix.
\end{proof}

Again, there is an analogue for operator-valued functions. The proof above requires only minimal changes.
\begin{theorem}\label{therealznthmop}
Let $\calh, \calk$ be Hilbert spaces.

If $\ph\in\ess(G; \calh,\calk)$ then there exist a Hilbert space $\calm$, a unitary operator $U$ on $\calm$ and a unitary operator
\beq\label{colligop}
\bbm A & B \\ C& D \ebm : \calh \oplus\calm \to \calk\oplus\calm
\eeq
such that, for all $s\in G$,
\beq\label{thisisitop}
\ph(s)=A+Bs_U(1-Ds_U)\inv C.
\eeq

Conversely, if a Hilbert space $\calm$, a contraction $T$ on $\calm$ and a contraction
\beq\label{colligbisop}
\bbm A & B \\ C & D \ebm : \calh\oplus\calm \to \calk\oplus\calm
\eeq
are given,  then the function $\ph:G\to \call(\calh,\calk)$ defined by
\beq\label{lfsTop}
\ph(s)=A+ Bs_T(1-Ds_T)\inv C
\eeq
belongs to $\ess(G;\calh,\calk)$.
\end{theorem}

\section{A second model formula for $G$ and spectral domains  }\label{secondmodel}
The model formula in Section \ref{model} has an alternative expression as an integral formula.

We shall need the rational functions 
\[
\Phi_\omega (s)=\frac{2\omega s^2- s^1 }{2-\omega s^1},\qquad s \in G,
\]
for $\omega\in\T$ (in the notation of the proof of Theorem \ref{modelGthm}, $\Phi_\omega(s)=f_s(\omega)$).  These functions have been used in many papers on $G$.  By \cite[Theorem 2.1]{AY04}, each $\Phi_\omega$ maps $G$ into $\D$.

Now invoke the spectral theorem to rewrite the model formula \eqref{modelform}.    Consider a function $\ph\in\ess(G)$.  By Theorem \ref{modelGthm}, $\ph$ has a $G$-model $(\calm, T, u)$ in which $T$ is a unitary operator on $\calm$.  By the spectral theorem,
\[
T = \int_\T \omega  \; \dd E(\omega),
\]
for some $\mathcal{L}(\m,\m)$-valued spectral measure $E$ on $\T$.   Thus, for $s\in G$,
\begin{align*}
s_T &= (2s^2T-s^1)(2-s^1T)\inv \\
	&=\int_\T \Phi_\omega(s) \; \dd E(\omega),
\end{align*}
and therefore
\[
1-t_T^* s_T = \int_\T 1- \overline{\Phi_\omega(t)}\Phi_\omega(s) \; \dd E(\omega).
\]
On combining this formula with Theorem \ref{modelGthm} we obtain the following statement.
\begin{theorem}\label{integralform}
 Let $\ph:G \to \C$ be a function. Then $\ph \in \ess (G)$ if and only if there exist a Hilbert space $\calm$, an $\mathcal{L}(\calm,\calm)$-valued spectral measure $E$ on $\T$ and an analytic map $u:G \to \m$ such that, for all $s,t\in G$,
\[
1-\overline{\ph(t)}\ph(s) = \int_{\T} \left( 1-\overline{\Phi_\omega(t)}\Phi_\omega(s)\right) \, \ip{\dd E(\omega)  u(s)}{u(t)}.
\]
\end{theorem}

One advantage of the integral form of the model formula is that it instantly yields a criterion for $G$  to be a spectral domain of a commuting pair of operators.  We recall the meaning of this notion.
\begin{definition} 
If $T$ is a $d$-tuple of pairwise commuting operators and $U$ is an open set in $\C^d$ we say that $U$ is a \emph{spectral domain for $T$} if $\sigma(T)\subset U$ and
\[
 \ph \in \ess(U) \implies \norm{\ph(T)} \le 1.       
\]
\end{definition}

The following statement is contained in \cite[Theorem 1.2]{AY}.
\begin{theorem} 
Let $S=(S^1,S^2)$ be a commuting pair of operators acting on a Hilbert space with $\sigma(S) \subset G$. Then $G$ is a spectral domain for $S$ if and only if
\[
 \norm{\Phi_\omega(S)} \le 1 \quad \mbox { for  all }\omega \in \T.
\]
\end{theorem}
\begin{proof} Since $\Phi_\omega \in \ess(G)$, the condition is obviously necessary.

Conversely, assume that $\norm{\Phi_\omega(S)} \le 1$ for all $\omega \in \T$. We need to show that $G$ is a spectral domain for $S$, i.e., that
\[
\norm{\ph(S)} \le 1
\]
whenever $\ph \in \ess(G)$.

But if $\ph \in \ess(G)$, it follows from Theorem \ref{integralform} that $1-\overline{\ph(t)}\ph(s)$ can be uniformly approximated by convex combinations of functions of the form
\[
\overline{f(t)}\Big(1-\overline{\Phi_\omega(t)}\Phi_\omega(s)\Big)f(s)
\]
where $\omega \in \T$ and $f$ is holomorphic on $G$. It follows that
$1-\ph(S)^*\ph(S)$ can be approximated in the operator norm by operators of the form
\[
f(S)^*\Big(1-\Phi_\omega(S)^*\Phi_\omega(S)\Big)f(S).
\]
Since these operators are positive, it follows that $1-\ph(S)^*\ph(S)$ is positive, that is, $\norm{\ph(S)} \le 1$.
\end{proof}

\section{A Pick theorem for $G$}\label{pickthm}
A standard application of realization formulae is to prove Pick-type theorems, which provide necessary and sufficient
conditions for the solvability of interpolation problems.  For example, the realization formula for the Schur class of the
bidisc in \cite{ag1990} yields the following criterion for analytic interpolation from $\D^2$ to $\D^-$.

Let $\la_1, \dots,\la_n$ be distinct points in $\D^2$ and let $w_1,\dots,w_n$ belong to $\D^-$.   There exists a function $\ph$ in $\ess(\D^2)$ such that $\ph(\la_j)=w_j$ for $j=1,\dots, n$ if and only if there exist positive semidefinite $n\times n$ matrices $a^1= [a^1_{ij}]_{i,j=1}^n$ and
$a^2=[a^2_{ij}]_{i,j=1}^n$ such that 
\be\label{bidiscpick}
1-\overline{w_i}w_j=  a^1_{ij}(1-\overline{\la_i^1}\la_j^1) + a^2_{ij} (1-\overline{\la_i^2}\la_j^2)
\ee
for $i,j=1,\dots,n$.

This result reduces the interpolation problem to the feasibility of a linear matrix equality for $n\times n$ matrices, a task which can be efficiently solved by standard engineering packages such as Matlab.

Consider the analogous problem in which the bidisc is replaced by the symmetrized bidisc.  Given distinct points $s_1, \dots,s_n$ in $G$
and target points $w_1,\dots,w_n$ in $\D^-$, we wish to determine whether there exists an analytic function $\ph: G\to \D^-$ such that
$\ph(s_j)=w_j$ for $j=1,\dots,n$.  One way to solve such an interpolation problem is to lift it to the bidisc.  Let $\mu_1,\dots,\mu_m$ be the preimages in $\D^2$ of the points $s_1,\dots,s_n$ under the natural map $\pi: \D^2\to G$ given by 
\[
\pi(\mu)=(\mu^1+\mu^2, \mu^1\mu^2).
\]
For any $s \in G$, the set $\pi\inv\{s\}$ comprises either one or two points, and therefore $n \leq m\leq 2n$.  It is easily seen that our
interpolation problem $s_j \mapsto w_j$ for $G$ is equivalent to the lifted problem $\mu_j \mapsto w_{j'}$ on $\D^2$, where $j'$ is chosen so 
that $1\leq j'\leq n$ and $\pi(\mu_j)=s_{j'}$.   The Pick criterion  \eqref{bidiscpick} applies to the lifted problem; since this criterion is necessarily symmetric with respect to the transposition map $\si$, it can be rewritten in terms of the symmetrized variables $s_j$ (and $w_j$).

However, the model in Theorem \ref{modelGthm} permits us to obtain directly a criterion for interpolation from $G$ to $\D^-$ in terms of the symmetrized variables.
\begin{theorem}\label{pickG}
Let $s_1, \dots,s_n$ be distinct points in $G$ and let $w_1,\dots,w_n\in \D^-$.   There exists an analytic function $\ph:G\to \D^-$ such that
$\ph(s_j)=w_j$ for $j=1,\dots,n$ if and only if there exist a Hilbert space $\calm$, a contraction $T$ on $\calm$ and vectors $v_1,\dots,v_n \in\calm$ such that 
\be\label{symvar}
1-\overline{w_i}w_j = \ip{\left(1-(s_i)_T^*(s_j)_T\right) v_j}{v_i}_\calm
\ee
for $i,j=1,\dots,n$.
\end{theorem}
\begin{proof}
{\em Necessity}.  Suppose an interpolating function $\ph\in\ess(G)$ exists.  By Theorem \ref{modelGthm}, $\ph$ has a $G$-model, that is, there exist a Hilbert space $\calm$, a contraction $T$ on $\calm$ and an analytic map $u:G\to\calm$ such that, for all $s,t\in G$,
\beq
 1-\overline{\ph(t)}\ph(s)= \ip{ (1-t_T^* s_T) u(s)}{u(t)}_\calm.
\eeq
On choosing $s=s_j$, $t=s_i$ and $v_i=u(s_i)$ for $i=1,\dots, n$ we deduce that equation \eqref{symvar} holds for all $i,j$.

{\em Sufficiency}.  Suppose that $\calm, T, v_1,\dots,v_n$ exist such that equation \eqref{symvar} holds for each $i,j$, as in the statement of the theorem.  Rearrange the equation to obtain
\[
1+\ip{(s_j)_Tv_j}{(s_i)_Tv_i}_\calm = \overline{w_i}w_j+\ip{v_j}{v_i}_\calm
\]
for all $i,j$.  This means that the family of vectors $(1,(s_j)_Tv_j), \; j=1,\dots,n$, in $\C\oplus\calm$ has the same Gramian as the family
$(w_j,v_j), \; j=1,\dots,n$, also in $\C\oplus\calm$.   Hence there exists an isometry
\[
L:\spa \{ (1,(s_j)_Tv_j) | j=1,\dots,n\} \to  \spa \{ (w_j,v_j) | j=1,\dots,n\}
\]
such that $L(1,(s_j)_Tv_j)= (w_j,v_j)$ for each $j$.  Extend $L$ to a contraction $L^\sharp$ mapping $\C\oplus \calm$ to itself.  Then $L^\sharp$ is expressible as a block operator matrix of the form
\[
L^\sharp \sim \bbm A & 1\otimes\beta \\ \ga\otimes 1& D\ebm
\]
for some $A\in\C$, vectors $\beta,\ga \in\calm$ and operator $D$ on $\calm$.
Since $L^\sharp (1,(s_j)_Tv_j)= (w_j,v_j)$ for each $j$,
\begin{align*}
A+\ip{(s_j)_Tv_j}{\beta}&=w_j,\\
\ga+ D(s_j)_Tv_j &= v_j.
\end{align*}
Thus
\[
v_j=(1-D(s_j)_T)\inv\ga
\]
and
\be\label{wj}
A+\ip{(s_j)_T(1-D(s_j)_T)\inv\ga}{\beta} = w_j
\ee
for each $j$.

Define a function $\ph:G\to\C$ by
\[
\ph(s)= A+\ip{s_T(1-Ds_T)\inv\ga}{\beta}.
\]
By Theorem \ref{modelGthm}, $\ph\in\ess(G)$, and by equation \eqref{wj},
\[
\ph(s_j)=w_j \quad\mbox{ for } j=1,\dots,n.
\]
\end{proof}
\begin{remark}   \rm
One can replace `there exists a contraction $T$' in the statement of Theorem \ref{pickG} by `there exists a unitary operator $T$'.
\end{remark}

Another criterion for the solvability of a finite interpolation problem in $\ess(G)$ is given in \cite[Theorem 6.1]{BS}.  It is shown that, in the situation of Theorem \ref{pickG}, a desired interpolating function $\ph\in\ess(G)$ exists if and only if there exists a $C(\D^-)^*$-valued positive semidefinite kernel on $\{s_1,\dots,s_n\}$ such that an analogue of equation \eqref{symvar} holds.

We conclude with an observation about the question raised in the introduction: is there a continuous {\em operator}-valued function $F$ on the closure $G^-$ of $G$, holomorphic in $G$, such that 
\beq
G=\{s\in G^-: \|F(s)\| < 1\}?
\eeq
Since a point $s\in\C^2$ belongs to $G$ if and only if $|f_s(\la)| < 1$ for all $\la\in\D$, one could try 
\[
F(s)= \diag_{n\geq 1}[f_s(\la_n)],
\]
where $(\la_n)$ is a dense sequence in $\D$.  It is then true that $F$ is well defined on $G^-$ and $G=\{s\in G^-: \|F(s)\| < 1\}$, but $F$ is discontinuous as a map from $G^-$ to the space of bounded linear operators on $\ell^2$ with the operator norm.  Indeed, for any $\omega\in\T$ and $0<r<1$,
\begin{align*}
\|F(2\bar\omega,\bar\omega^2)-F(2r\bar\omega,r\bar\omega^2)\|&=(1-r)\sup_n\left|\frac{\la_n}{1-r\la_n\bar\omega}\right| \\
	&=1.
\end{align*}
Thus $F$ is discontinuous at every point $(2\bar\omega,\bar\omega^2) \in G^-$ for $\omega\in\T$.
Indeed $F$ is even discontinuous at these points with respect to the weak operator topology on the space of bounded linear operators on $\ell^2$.
We leave open the question of whether there exists a continuous holomorphic operator-valued defining function for $G$.

\end{document}